\documentclass[a4paper]{article}
\usepackage[colorinlistoftodos,textsize=small]{todonotes}
\usepackage[margin=3.2cm]{geometry}
\usepackage{amsmath,amssymb}
\usepackage{mathptmx}
\usepackage{helvet}
\usepackage{courier}
\usepackage{makeidx}
\usepackage{graphicx}
\usepackage{multicol}
\usepackage[bottom]{footmisc}
\usepackage{cite}
\usepackage{fancyvrb}
\usepackage{subcaption}
\usepackage{hyperref}
\usepackage{pgfplots}

\usepackage{tikz}
\usetikzlibrary{decorations,arrows,calc,shapes,arrows,positioning,shadows.blur,shapes.symbols}
\usepackage[siunitx,europeanresistors]{circuitikz}
\usetikzlibrary{spy,backgrounds,decorations,arrows,calc,shapes,arrows,positioning,shadows.blur,shapes.symbols}

\newenvironment{Proof}{\noindent{\scshape{\textbf{Proof:}}}$ $\ }

\newtheorem{Defn}{Definition}
\newtheorem{Thm}[Defn]{Theorem}

\newtheorem{Ass}[Defn]{Assumption}
\newtheorem{Lem}[Defn]{Lemma}
\newtheorem{Rem}[Defn]{Remark}

\newcommand{\R}{\mathbb{R}}

\begin{document}

\title{Waveform relaxation for\\low frequency coupled field/circuit\\differential-algebraic models of index~2}
\author{
  Idoia Cortes Garcia$^{1}$
  \and
  Jonas Pade$^{2}$
  \and
  Sebastian Sch\"ops$^{1}$
  \and
  Caren Tischendorf$^{2}$\\[1em]
  $1$~Technical University of Darmstadt, CEM Group,\\ 
  Schlossgartenstrasse 8, 64289 Darmstadt, Germany\\[1em]
  $2$~Department of Mathematics, Humboldt University of Berlin,\\
  Rudower Chaussee 25, 12489, Berlin, Germany%
}

\maketitle

\abstract{Motivated by the task to design quench protection systems for superconducting magnets in particle accelerators we address a coupled field/circuit simulation based on a magneto-quasistatic field modeling. We investigate how a waveform relaxation of Gau\ss{}-Seidel type performs for a coupled simulation when circuit solving packages are used that describe the circuit by the modified nodal analysis. We present sufficient convergence criteria for the coupled simulation of FEM discretised field models and circuit models formed by a differential-algebraic equation (DAE) system of index 2. In particular, we demonstrate by a simple benchmark system the drastic influence of the circuit topology on the convergence behavior of the coupled simulation.}

\section{Introduction}
\label{pade:sec_field/circuit}
Lumped circuit models, such as modified nodal analysis (MNA), are well-established in electrical engineering. 
However, they neglect the spatial dimension and therefore distributed phenomena like the skin effect. For certain devices, this may lead to inaccuracies of unacceptable magnitude in the simulation, e.g. for electric machines \cite{Salon_1995aa}
or the quench protection system of superconducting magnets in particle accelerators \cite{Bortot_2018ab}.
These cases call for field/circuit coupling \cite{Schops_2010aa}, \cite{Cortes-Garcia_2017ab}.
To solve such coupled systems, it is often advisable to use waveform relaxation (WR) \cite{Lelarasmee_1982ab}, since this iterative method allows for 
dedicated step sizes and suitable solvers for the different subsystems, and even for the use of proprietary blackbox solvers.
The coupled field/circuit model considered here is a DAE in the time domain after space discretisation of the field system. It is well-known 
that WR can suffer from instabilities for DAEs unless an additional contraction criterion is satisfied \cite{Lelarasmee_1982ab}, \cite{Pade_2018aa}.
This work presents coupled field/circuit models, which are DAEs of index $2$ \cite{Estevez-Schwarz_2000aa}, for the case where WR is convergent and the case where it diverges.
Furthermore, generalizing a convergence criterion of \cite{Pade_2018aa}, a topological and easy-to-check criterion is provided. 
Finally, we present numerical simulations verifying the topological convergence criterion.

\section{Field/Circuit Model}
\label{pade:sec_field/circuit}
To describe the electromagnetic (EM) field part, we consider a magnetoquasistatic approximation of Max\-well's equations in a reduced magnetic vector potential formulation \cite{Emson_1988aa}.
This leads to the curl-curl eddy current partial differential equation (PDE).
The circuit side is formulated with the MNA \cite{Ho_1975aa}. For the numerical simulation of the coupled system, the method of lines is used with a finite element (FE) discretisation.
Altogether, this leads to a time-dependent coupled system of DAE initial value problems (IVPs), described by
\begin{align}
 \label{eq:em1} M\dot a+K(a)a-Xi_m&=0,& X^\top \dot a&=v_c,\\
 \label{eq:circ1} E(x)\dot x+f(t,x)&=Pi_m, & P^\top x-v_c&=0.
\end{align}
The first Equation \eqref{eq:em1} represents the space-discrete field model based on the matrices
\begin{align}\label{eq:fem_matrices}
 (M)_{ij} = \int_{\Omega} \sigma \omega_i \cdot \omega_j\,\mathrm{d}\mathrm{V}, &&
 (K(a))_{ij} = \int_{\Omega} \nu(a) \nabla\times\omega_i \cdot   \nabla\times \omega_j\,\mathrm{d}\mathrm{V}\;,
\end{align}
which follow from the Ritz-Galerkin approach using a finite set of N\'ed\'elec basis functions $\omega_i$ \cite{Monk_2003aa} defined on the domain $\Omega$; $\sigma$ denotes the space-dependent electric conductivity and $\nu(a)$ the magnetic reluctivity that can additionally depend nonlinearly on the unknown magnetic vector potential $a$.
The excitation matrix is computed from a winding density function ${\chi}_j$ modelling the $j$-th
stranded conductor \cite{Schops_2011ac} as
\begin{align}
	 (X)_{ij} = \int_{\Omega} {\chi}_j\cdot \omega_i \,\mathrm{d}\mathrm{V}\;.
\end{align}
\begin{Defn}
 A function $f:\R^n\to\R^n$ is \emph{strongly monotone} and a square matrix $M(x)$ is \emph{uniformly positive definite}, if
 \begin{align*}
  \exists \mu_f:&& (x_2-x_1)^\top(f(x_2)-f(x_1))\geq &\mu_f\|x_2-x_1\|^2,
    \ &\forall x_1,x_2\in\mathbb R^n,\\
    \exists \mu_M:&& y^\top M(x)y \geq &\mu_M\|y\|^2,\ &\forall x\in\R^n,y\in\mathbb R^m.
 \end{align*}
\end{Defn}
The space-discretization is supposed to meet the following properties.
\begin{Ass}
 \label{as:ind1}
 It holds (a) $M$ is symmetric, (b) the matrix pencil $\lambda M+K$ is symmetric and positive definite for $\lambda>0$, (c) $X$ has full column rank and (d) 
 the function $a\mapsto K(a)a$ is strongly monotone.
\end{Ass}
The assumptions are in agreement with previous formulation in the literature, e.g. \cite{Cortes-Garcia_2020ab, Schops_2011ac}. The first Assumption~\ref{as:ind1}a) follows naturally if a Ritz-Galerkin formulation \eqref{eq:fem_matrices} is chosen. The second Assumption~\ref{as:ind1}b) will be guaranteed by appropriate boundary and gauging conditions. Thirdly, the full column rank Assumption~\ref{as:ind1}c) follows from the fact that the columns are discretisations of different coils that are located in spatially disjoint subdomains.
Finally, the monotonicity Assumption ~\ref{as:ind1}(d) follows from the strong monotonicity of the underlying nonlinear material law, i.e. the BH-curve \cite{Pechstein_2004aa}.
In general, the field model is a multiport element such that the circuit coupling is established via multiple currents and voltages, i.e., vector-valued $i_m$ and $v_c$. However, for simplicity of notation we assume a two-terminal device in the following.

\medskip

The circuit Equation \eqref{eq:circ1} can be expanded into
\begin{align}
  \label{eq:MNA}\small{
  E(x) =\begin{pmatrix}\mathcal L_C(e)&0&0\\0&-L(i_L)&0\\ 0&0&0\end{pmatrix},\;
  f(t,x) =  \begin{pmatrix}g_R(e)+A_L i_L+A_V i_V+q_i(t)\\A_L^\top e\\A_V^\top e - q_v(t)\end{pmatrix},\;
  P=\begin{pmatrix}A_m\\0\\0\end{pmatrix}}
\end{align}
using the definitions $\mathcal L_C(e):=A_C C(A_C^\top e)A_C^\top$, $g_R(e):=A_R g(A_R^\top e)$ and $x=(e,i_L,i_V)$ where $A_\star$ are the usual incidence matrices and $L(\cdot)$,
$C(\cdot)$ 
are state-dependent square matrices describing inductances and capacitances. The function $g(\cdot)$
describes the voltage-current relation of resistive elements. Finally, $x$ collects all node potentials $e$, currents through branches with voltage sources $i_V$ and inductors $i_L$. 
The circuit system shall fulfill the following properties:
\begin{Ass}
 \label{as:circas}
  It holds (a) $g$, $C$ and $L$ are Lipschitz continuous, 
    $g$ is strongly monotone and $C,\ L$ are uniformly positive definite, (b) $q_i$ and $q_v$ are continuously differentiable, (c) $A_V$ has full column rank and $\begin{pmatrix}A_C&A_V&A_R&A_L\end{pmatrix}$ has
    full row rank.
\end{Ass}
Assumption ~\ref{as:circas}a) reflects the global passivity of the respective elements \cite{Matthes_2012aa}. 
Considering well-known relations between incidence matrices and circuit topology, Assumption ~\ref{as:circas}c) 
excludes the electrically forbidden configurations of loops of voltage sources and cutsets of current sources \cite{Estevez-Schwarz_2000aa}.
\section{Waveform Relaxation and Convergence}
 \label{Pade:sec_WR}
We consider the Gau\ss-Seidel WR method. Applied to the coupled system \eqref{eq:em1}-\eqref{eq:circ1}, this yields the scheme
\begin{align}
\label{eq:WR1}M\dot a^k+K(a^k)a^k-Xi_m^k&=0,&\ X^\top \dot a^k&=v_c^{k-1},\\
\label{eq:WR4}E(x^k)\dot x^k+f(t,x^k)&=Pi_m^k,&\ P^\top x^k-v_c^k&=0.
\end{align}
The superscript $k$ denotes the iteration index. A common choice for the inital guess $v_c^{0}$ is constant extrapolation of the initial value.

We shall proceed as follows:\medskip

\begin{enumerate}
 \item Lemmata \ref{lem:EMdec} and \ref{lem:MNAdec} provide a DAE-decoupling of the EM field DAE \eqref{eq:em1} and the 
 MNA DAE \eqref{eq:circ1}, respectively.
 \item Definition \ref{def:CVR} introduces the concept of parallel CVR paths.
 Assuming their existence and exploiting the previous decoupling Lemmata, Lemma \ref{lem:inhODE} yields a DAE-decoupling of the coupled WR iteration 
 \eqref{eq:WR1}-\eqref{eq:WR4}. Notably, it reveals the structure of its inherent ODE, given by $\phi$ in Equation \eqref{eq:WRdec}.
 \item The convergence Theorem \ref{thm:crit} is a simple consequence of the previous Lemmata; it shows that 
 the existence of parallel CVR paths guarantees convergence of the WR scheme \eqref{eq:WR1}-\eqref{eq:WR4}.
\end{enumerate}
For visual reasons, we shall write column vectors as $(a,b,c)$.
\begin{Lem}\label{lem:EMdec}
Let Assumption \ref{as:ind1} hold. Then, for a given source term $v_c$, there exists a coordinate transformation 
 $(w,u)=T^{-1}a$
 and a system of the form
 \begin{align}
  \label{eq:equiv}\dot u+A_1u&=A_2v_c,&\ w&=Bu,&\ i_m=G_1u+G_2v_c
 \end{align}
 such that $(a,i_m)$ solves Equation \eqref{eq:em1} if and only if $(u,w,i_m)$ solves Equation \eqref{eq:equiv}.
\end{Lem}
\begin{Proof}
For better readability and shortness we present the proof only for the slightly more restrictive case where $X^\top M=0$, which is usually satisfied.

 We equivalently transform the field DAE with new coordinates 
$T\alpha=a$:
\begin{align}
\begin{split}\label{eq:trafoDAE}
T^\top M T\dot\alpha+ T^\top K(T\alpha)T\alpha- T^\top Xi_m&=0,\\
 X^\top T\dot\alpha&=v_c.
\end{split}
 \end{align}
The transformation matrix $T:=(T_{\ker}\ X\ T_\perp)$ is constructed such that the columns of $T_{\ker}$ and $T_\perp$ form a basis of $\ker M\cap\ker X^\top$ and $(\ker M)^\perp$,
respectively. It is nonsingular indeed, since its construction and Assumption ~\ref{as:ind1} combined with $X^TM=0$ guarantee that
$\text{im}X\perp\text{im}T_{\ker}$ and $\text{im}T_\perp\perp\text{im}(T_{\ker}\ X)$.

With $\alpha=(w,u)$ and $u=(u_1,u_2)$, the transformed DAE \eqref{eq:trafoDAE} has the detailed form
\begin{align*}
 \underline{T_{\ker}^\top K(T\alpha)T_{\ker}}w+T_{\ker}^\top K(T\alpha)(X\ \ T_\perp)u&=0,\\
 X^\top K(T\alpha)T\alpha-\underline{X^\top X}i_m&=0,\\
  \underline{T_\perp^\top M T_\perp} \dot u_2+T_\perp^\top K(T\alpha)T\alpha&=0,\\
  \underline{X^\top X}\dot u_1&= v_c.
\end{align*} 
The underlined matrices are nonsingular due to Assumption ~\ref{as:ind1}, and Equation \eqref{eq:equiv} is obtained by inversion and insertion.
\end{Proof}
\begin{Defn}\label{def:CVR}
 A \emph{CVR path} in a circuit is a path which consists of only capacitances, voltages sources and resistances.
 An element has a \emph{parallel CVR path}, if its incident nodes are connnected by a CVR path.
\end{Defn}
\begin{Lem}\label{lem:MNAdec}
Let Assumption ~\ref{as:circas} hold. Then, for a given source term $i_m$, there exists
a coordinate transformation
 $(y,z_1,z_2)=T^{-1}x$
 and a system of the form
 \begin{subequations}\label{eq:MNAdec}
 \begin{align}\label{eq:MNAdec1}
  \dot {y}&= f_0(t,y,z,z_2,u),&
    z_1&=g_1(t,y,z_2,\dot z_2,u),&
    z_2&=g_2(t)+QPi_m,\\
    &&&&v_c&=P^\top T(y,z_1,z_2)
 \end{align}\end{subequations}
 with $f_0,g_1,g_2$ uniformly globally Lipschitz continuous $\forall t$ and $g_2\in C^1$ such that
 \begin{enumerate}
  \item $(x,v_c)$ solves Equation \eqref{eq:circ1} if and only if $(y,z_1,z_2,v_c)$ solves Equation \eqref{eq:MNAdec},
  \item $QP=0$ if each EM field element has a parallel CVR-path.
 \end{enumerate}
 \end{Lem}
A detailed proof can be found in \cite{Pade_2020aa},
where $Q$ is shown to have the form $(Q_1\ *\ *)$ with $\text{im}Q_1=\ker (A_C\ A_V\ A_R)^\top$.
Hence, if each field element has a parallel CVR-path, each column of $A_m$ 
can be written as a sum of columns of $(A_C\ A_V\ A_R)$ and it follows $Q_1A_m=0$, thus $QP=0$.
\begin{Lem}\label{lem:inhODE}
Let Assumptions \ref{as:ind1} and \ref{as:circas} hold. 
If each EM field element has a parallel CVR path, then
there exists a coordinate transformation
 $(r,s)=T^{-1}(a,x)$
 and a system of the form
 \begin{align}\label{eq:WRdec}
   \dot s^k=\phi(t,s^k,s^{k-1}),\qquad 
    r^k=\varphi(t,s^k)
 \end{align}
 with $\phi$ uniformly globally Lipschitz continuous $\forall t$ and $\phi,\varphi$ continuous
 such that $(a^k,i_m^k,x^k,v_c^k)$ solves Equations \eqref{eq:WR1}-\eqref{eq:WR4} if and only if $(s^k,r^k)$ solves Equation \eqref{eq:WRdec}. 
\end{Lem}
\begin{Proof}
We apply Lemmata ~\ref{lem:EMdec},~\ref{lem:MNAdec} to the iterated subsystems \eqref{eq:WR1},\eqref{eq:WR4}. This yields
an equivalent system
\begin{align}\label{eq:aux1}
  \dot u^k&=-A_1u^k+A_2v_c^{k-1}, &
  w^k&=Bu^k, &
  i_m^k&=G_1u^k+G_2v_c^{k-1},\\
  \label{eq:aux2}
  \dot y^k &= f_0(t,y^k,z^k,z_2^k,u^k), &
  z_1^k&=g_1(t,y^k,z_2^k,\dot z_2^k,u^k), &
  z_2^k&=g_2(t),\\
  && &&v_c^k &=P^\top T(y^k,z_1^k,z_2^k). \label{eq:aux3}
    \end{align}
Since each field element has a parallel CVR path, $z_2^k=g(t)$ does not depend on $u^k$ anymore.

We insert $v_c^{k-1}=P^\top x^{k-1}=P^\top T(y^{k-1},z_1^{k-1},z_2^{k-1})$
and $z_1^{k-1}$ and $z_2^{k-1}$ therein to obtain, 
with $\tilde g_1(t,y^{k-1},u^{k-1})=g_1(t,y^{k-1},g_2(t),\dot g_2(t),u^{k-1})$,
\begin{align*}
 \dot u^k=\phi_2(t,u^k,y^k,u^{k-1},y^{k-1})\quad :=-A_1u^k+A_2P^\top T(y^{k-1},\tilde g_1(t,y^{k-1},u^{k-1}),g_2(t)).
\end{align*}
Insertion of $z_1^k,z_2^k,\dot z_2^k$ into $f_0$ yields 
\begin{align*}
 \dot y^k=\phi_1(t,u^k,y^k)\quad :=f_0(t,y^k,g_1(t,y^k,g_2(t),\dot g_2(t),u^k),g_2(t),u^k).
\end{align*}
Hence, defining $s^k:=(u^k,y^k)$ and $\phi:=(\phi_1,\phi_2)$, the sequence $(u^k,y^k)$ is given implicitly by an ODE recursion of the form
 $\dot s^k=\phi(t,s^k,s^{k-1})$.
 
 The algebraic constraint of Equation \eqref{eq:WRdec} is obtained 
 with $r^k=(w^k,i^k,z_1^k,z_2^k,v_c^k)$, $s^k=(u^k,y^k)$ and 
 \begin{align*}
  \varphi(t,s)=(Bu,Gu,g_1(t,y,g_2(t),\dot g_2(t),u),g_2(t)).
 \end{align*}
Clearly, $(s^k,r^k)$ solves Equation \eqref{eq:WRdec} if and only if $\tilde\alpha^k:=(u^k,w^k,i_m^k,y^k,z_1^k,z_2^k,v_c^k)$ solves 
Equations \eqref{eq:aux1}-\eqref{eq:aux3}, and $\tilde\alpha^k$ solves \eqref{eq:aux1}-\eqref{eq:aux3} if and only if $(a^k,i_m^k,x^k,v_c^k)$
solves Equations \eqref{eq:WR1}-\eqref{eq:WR4}.
 \end{Proof}
We deduce the main result of this work:

\begin{Thm}\label{thm:crit}
 If each EM field element of the coupled system \eqref{eq:em1}-\eqref{eq:circ1} has a parallel CVR path, then the WR scheme \eqref{eq:WR1}-\eqref{eq:WR4} 
 is uniformly convergent to the
 exact solution of \eqref{eq:em1}-\eqref{eq:circ1}. 
\end{Thm}\label{thm:conv}
\begin{Proof}
The ODE part of Equation \eqref{eq:WRdec} is a WR scheme for ODEs with Lipschitz continuous vector field $\phi$.
It is well-known that such schemes are unconditionally convergent on bounded time intervals \cite{Lelarasmee_1982ab}.
The convergence of $s^k$ clearly implies the convergence of $(s^k,r^k)$ defined by \eqref{eq:WRdec}.
Due to the equivalence provided by Lemma \ref{lem:inhODE}, 
it follows that the original scheme \eqref{eq:WR1}-\eqref{eq:WR4} is convergent. 
\end{Proof}
\begin{Rem}
The convergence result holds for arbitrary continuous initial guesses $x^0$ 
and for bounded intervals of arbitrary size, see e.g. \cite{Lelarasmee_1982ab},\cite{Pade_2020aa}.
\end{Rem}
\begin{figure}[t]
  \centering
  \vspace{-1em}
  \begin{tikzpicture}
  \node at (0,0) {\includegraphics[width=0.5\textwidth]{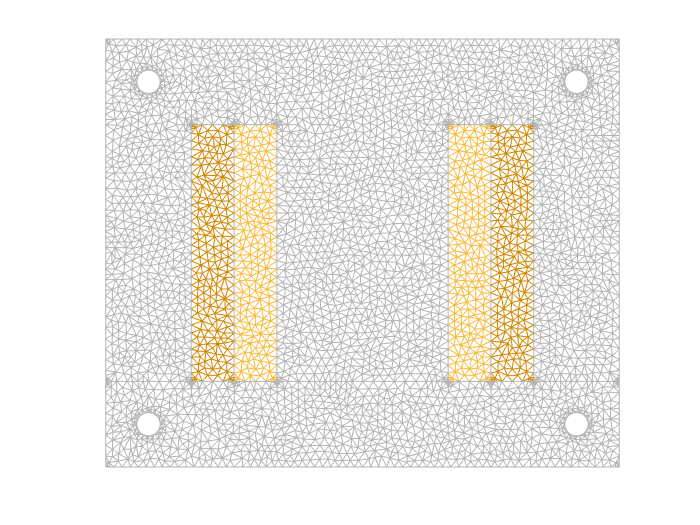}};
  \draw[-latex, line width = 1pt] (2.8,1.6) -- (1,0);
  \node[right] at (2.8,1.6) {Primary coil};
  \end{tikzpicture}
  \vspace{-2.5em}
  \caption{
    Single phase isolation transformer (`MyTransformer'), see
\cite{Meeker_2018aa}.
  }
  \label{fig:transformer_wr_index2}
\end{figure}
\begin{Rem} The MNA decoupling given in Lemma ~\ref{lem:MNAdec} shows that $g_1$ depends on $z_2$ and the derivative $\dot z_2$.
Hence, the system is most sensitive to perturbations of $z_2$.
The input of the EM field subsystem in the WR scheme is in fact a perturbation.
Therefore, the condition $QP=0$ from Lemma ~\ref{lem:MNAdec} is crucial to derive Theorem ~\ref{thm:conv}.
If at least one EM field element has no parallel CVR path, then $QP\neq 0$. Then, analogously to Lemma ~\ref{lem:inhODE} and its proof, we find
 $\dot s^k=\phi(t,s^k,s^{k-1},\dot s^{k-1})$,
which is guaranteed to converge only if $\phi$ is contractive in $\dot s^{k-1}$, see \cite{Lelarasmee_1982ab},\cite{Pade_2020aa}.
\end{Rem}
 \section{Numerical Examples}
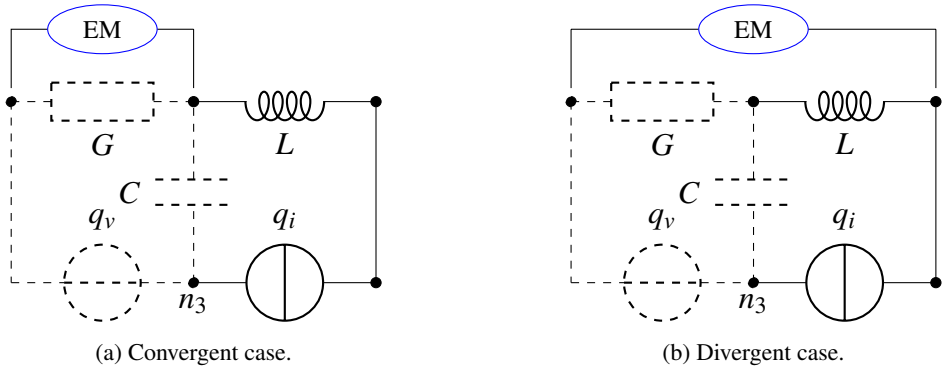
\begin{figure}[t]
  \begin{subfigure}[c]{0.5\textwidth}
    \centering
    \begin{circuitikz}[scale=1.2, transform shape]
    \draw[dashed] (0,0) to [V,-*,l=$q_v$] (2,0);
    \draw (2,0) to [I,-*,l=$q_i$] (4,0)
    to (4,2)
    to [L,*-*,l=$L$] (2,2);
    \draw[dashed] (2,2) to [R,-*,l=$G$] (0,2)
    to (0,0)
    (2,0) to [C,l=$C$] (2,2)
    ;
       \draw
      (0,2.15) to (0,2.8)
      to (0.4,2.8)
      (2,2.8) to (2,2.15)
      (2,2.8) to (1.6,2.8)
      ;
      \draw (2,0) node[anchor=north] {$n_3$};
      \draw[blue,fill=white] (1,2.8) ellipse (0.6cm and 0.27cm);
     \node at (1,2.8) {\footnotesize{EM}};
          \end{circuitikz}
    \caption{Convergent case.}
    \label{fig:convcirc}
  \end{subfigure}
  \begin{subfigure}[c]{0.5\textwidth}
     \centering
    \begin{circuitikz}[scale=1.2, transform shape]
    \draw[dashed] (0,0) to [V,-*,l=$q_v$] (2,0);
    \draw (2,0) to [I,-*,l=$q_i$] (4,0)
     to (4,2)
     to [L,*-*,l=$L$] (2,2);
    \draw[dashed] (2,2) to [R,-*,l=$G$] (0,2)
     to (0,0)
     (2,0) to [C,l=$C$] (2,2)
      ;
       \draw
      (0,2.15) to (0,2.8)
      (0,2.8) to (1.4,2.8)
      (4,2.8) to (4,2.15)
      (4,2.8) to (1.6,2.8)
      ;
      \draw (2,0) node[anchor=north] {$n_3$};
      \draw[blue,fill=white] (2,2.8) ellipse (0.6cm and 0.27cm);
     \node at (2,2.8) {\footnotesize{EM}};
          \end{circuitikz}
    \caption{Divergent case.}
    \label{fig:divcirc}
  \end{subfigure}
  \caption{Field/circuit coupling with model from Fig.~\ref{fig:transformer_wr_index2} (CVR path is dashed).\label{fig:circ}}
\end{figure}

 To illustrate the convergence behaviour of the WR scheme according to the derived criteria, we consider the toy
 example circuits in Figures~\ref{fig:convcirc}~and~\ref{fig:divcirc}. Both are described with MNA \eqref{eq:circ1} 
 and the (arbitrary) parameters $R=1\Omega$, $L=5$H, $C=1$F, $i_{\mathrm{s}}(t)=\sin(2t)+5\sin(20t)$ and $v_{\mathrm{s}}(t)=\sin(t)+\sin(20t)$
 are set. The eddy current Equation \eqref{eq:em1} is solved on the single phase isolation transformer shown in Figure~\ref{fig:transformer_wr_index2}.
 For simplicity, a zero current is imposed on the secondary coil (dark orange) and only the primary coil is coupled to the circuit.

The WR algorithm is applied on the simulation time window $\mathcal{I} = [0\;0.8]\,$s
and the internal time integration is
performed with the implicit Euler scheme with time step size $\delta t = 10^{-2}\,$s.

\begin{figure}
  \begin{subfigure}[c]{0.5\textwidth}
    \centering
    \begin{tikzpicture}
     \begin{axis}
     [width=0.9\linewidth,height=5cm,xlabel={{time t / s}},ylabel={{Potential $e$} / V},
     legend style={at={(0.1,1.1)},anchor=north west, ymax=250, ymin=-90, xlabel near ticks, ylabel near ticks}
     ]
     \addplot[mark=none, color=black, line width = 1pt] table [x=t, y=mon, col sep=comma] {convergent.csv};
     \addplot[dashed, color=blue, line width = 1pt] table [x=t, y=1, col sep=comma] {convergent.csv};
     \addplot[dashed, color=red, line width = 1pt] table [x=t, y=2, col sep=comma] {convergent.csv};
     \legend{\footnotesize{mon},\footnotesize{$k=1$},\footnotesize{$k=2$}}
     \end{axis}
    \end{tikzpicture}
    \caption{Convergent case.}
    \label{fig:convplot}
  \end{subfigure}
  \begin{subfigure}[c]{0.5\textwidth}
    \centering
    \begin{tikzpicture}
    \begin{axis}
    [width=0.9\linewidth,height=5cm,xlabel={{time t / s}},ylabel={{Potential $e$} / V},
     legend style={at={(0.1,1.1)},anchor=north west, ymax=250, ymin=-90, xlabel near ticks, ylabel near ticks}
     ]
    \addplot[mark=none, black, line width =1pt] table [x=t, y=mon, col sep=comma] {divergent.csv};
    \addplot[dashed, color=blue, line width = 1pt] table [x=t, y=1, col sep=comma] {divergent.csv};
    \addplot[dashed, color=red, line width = 1pt] table [x=t, y=2, col sep=comma] {divergent.csv};
    \legend{\footnotesize{mon},\footnotesize{$k=1$},\footnotesize{$k=2$}}
    \end{axis}
    \end{tikzpicture}
    \caption{Divergent case.}
    \label{fig:divplot}
  \end{subfigure}
  \caption{Monolithic ("mon") and WR solution for $k=1,2$ iterations.}
\end{figure}
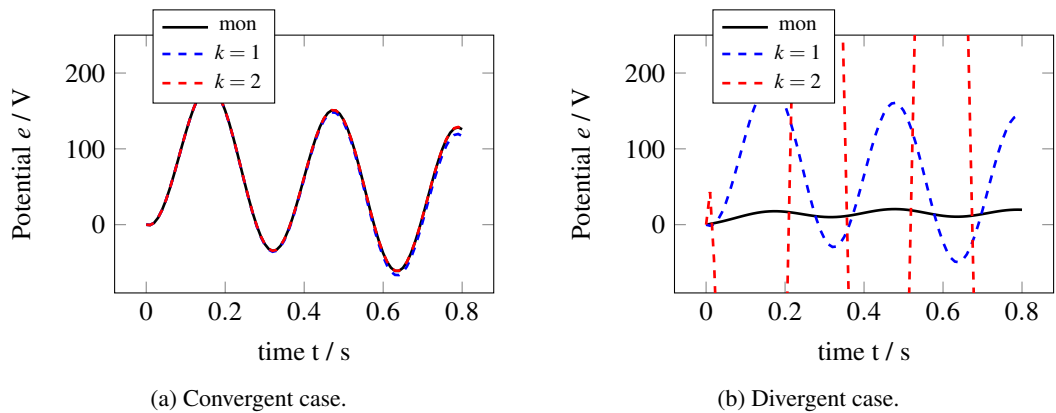

The theoretical result is illustrated by the successful simulation, see Figure~\ref{fig:convplot}, of the model shown in Figure~\ref{fig:convcirc} which satisfies the convergence criterion of Theorem~\ref{thm:crit}.
However, numerical simulations of the model shown in Figure~\ref{fig:divcirc} show that WR can diverge indeed if the criterion is not satisfied.

\section{Conclusions}
In this work, we have presented a space-discretised coupled field/circuit model, which is a DAE of index $2$, and a simulation of this model by means of WR.
Furthermore, we have provided an easy-to-check topological convergence criterion for a class of coupled DAE/DAE systems of index $2$.

\section*{acknowledgement}
This work is supported by the `Excellence Initiative' of the German Federal and State Governments, the Graduate School of CE at TU Darmstadt and DFG grant SCHO1562/1-2.
Further, we acknowledge financial support under BMWi grant 0324019E and by
DFG under Germany's Excellence Strategy – The Berlin Mathematics Research Center MATH+ (EXC-2046/1, ID 390685689).


\begin{thebibliography}{10}

\bibitem{Bortot_2018ab}
L. Bortot et. al.
\newblock {STEAM}: A hierarchical co-simulation framework for superconducting
  accelerator magnet circuits.
\newblock 28(3), 2018.

\bibitem{Cortes-Garcia_2020ab}
I. Cortes~Garcia, H. De~Gersem, and S. Sch\"{o}ps.
\newblock A structural analysis of field/circuit coupled problems based on a
  generalised circuit element.
\newblock 83(1):373--394, 2020.

\bibitem{Cortes-Garcia_2017ab}
I. Cortes~Garcia et. al.
\newblock Optimized field/circuit coupling for the simulation of quenches in
  superconducting magnets.
\newblock 2(1):97--104, 2017.

\bibitem{Emson_1988aa}
C.~R.~I. Emson and C.~W. Trowbridge.
\newblock Transient 3d eddy currents using modified magnetic vector potentials
  and magnetic scalar potentials.
\newblock 24(1):86--89, 1988.

\bibitem{Estevez-Schwarz_2000aa}
D: Est\'{e}vez~Schwarz and C. Tischendorf.
\newblock Structural analysis of electric circuits and consequences for {MNA}.
\newblock 28(2):131--162, 2000.

\bibitem{Ho_1975aa}
C.-W. Ho, A.~E. Ruehli, and P.~A. Brennan.
\newblock The modified nodal approach to network analysis.
\newblock 22(6):504--509, 1975.

\bibitem{Lelarasmee_1982ab}
E. Lelarasmee, A.~E. Ruehli, and A.~L. Sangiovanni-Vincentelli.
\newblock The waveform relaxation method for time-domain analysis of large
  scale integrated circuits.
\newblock 1(3):131--145, 1982.

\bibitem{Matthes_2012aa}
M. Matthes.
  Numerical Analysis of Nonlinear Partial Differential-Algebraic Equations: A Coupled and an Abstract Systems Approach.
  Logos Verlag Berlin GmbH, 2012.

\bibitem{Meeker_2018aa}
D. Meeker.
\newblock {\em Finite Element Method Magnetics}, version 4.2 (25feb2018 build)
  edition, 2018.
\newblock User's Manual.

\bibitem{Monk_2003aa}
P. Monk.
\newblock {\em Finite Element Methods for {Maxwell}'s Equations}.
\newblock Oxford University Press, 2003.

\bibitem{Pade_2018aa}
J. Pade and C. Tischendorf.
\newblock Waveform relaxation: a convergence criterion for
  differential-algebraic equations.
\newblock 81(4): 1327-1342, 2019.

\bibitem{Pade_2020aa}
J. Pade.
\newblock Convergence criteria for waveform relaxation on 
differential-algebraic systems: a topological approach for circuits. 
\newblock Ph.D. Thesis,
HU Berlin.
\newblock In Preparation, 2020

\bibitem{Pechstein_2004aa}
C. Pechstein.
  \newblock Multigrid-{Newton}-Methods For Nonlinear-Magnetostatic Problems,
\newblock Master's Thesis, University of Linz, 2004.

\bibitem{Salon_1995aa}
S.~J. Salon.
\newblock {\em Finite Element Analysis of Electrical Machines}.
\newblock Kluwer, 1995.

\bibitem{Schops_2011ac}
S. Sch\"{o}ps.
\newblock {\em Multiscale Modeling and Multirate Time-Integration of
  Field/Circuit Coupled Problems}.
\newblock {VDI} Verlag. Fortschritt-Berichte {VDI}, Reihe 21.

\bibitem{Schops_2010aa}
S. Sch\"{o}ps, H. De~Gersem, and A. Bartel.
\newblock A cosimulation framework for multirate time-integration of
  field/circuit coupled problems.
\newblock 46(8):3233--3236, 2010.

\end{thebibliography}
\end{document}